\documentclass[12pt]{article}
\usepackage{amscd, amssymb, latexsym, amsmath}

\newtheorem{thm}{Theorem}
\newtheorem{cor}{Corollary}
\newtheorem{pro}{Proposition}
\newtheorem{lem}{Lemma}
\newtheorem{dfn}{Definition}

\newenvironment{proof}
{\noindent {\em Proof.}}
{\hfill $\Box$}

\numberwithin{thm}{section}
\numberwithin{cor}{section}
\numberwithin{pro}{section}
\numberwithin{lem}{section}
\numberwithin{dfn}{section}
\numberwithin{rem}{section}
\numberwithin{equation}{section}

\newcommand{\R}{\mathbb R}

\newcommand{\h}{\mathbb H}

\newcommand{\ct}{\check{t}}
\newcommand{\cx}{\check{x}}
\begin{document}
\title{Deforming Area Preserving Diffeomorphism
of Surfaces by Mean Curvature Flow}
\author{Mu-Tao Wang}
\date{September 21, 2000, this version April 15, 2001}
\maketitle

\begin{abstract}
Let $f:\Sigma_1\rightarrow \Sigma_2 $ be an area preserving
diffeomorphism between compact Riemann surfaces of constant
curvature. The graph of $f$ can be viewed as a Lagrangian
submanifold in $\Sigma_1\times \Sigma_2$. This article discusses a
canonical way to deform $f$ along area preserving diffeomorphisms.
This deformation process is realized through the mean curvature
flow of the graph of $f$ in $\Sigma_1\times \Sigma_2$. It is
proved that the flow exists for all time and the map converges to
a canonical map. In particular, this gives a new proof of the
classical topological results that $O(3)$ is a deformation retract
of the diffeomorphism group of $S^2$ and the mapping class group
of a Riemman surface of positive genus is a deformation retract of
the diffeomorphism group .
\end{abstract}

\section{Introduction}

The mean curvature flow is an evolution process under which a
submanifold evolves in the direction of its mean curvature vector.
It represents the most effective way to decrease the volume of a
submanifold. The codimension one case of the mean curvature flow
has been studied extensively while very little is known in the
higher codimension case. The multidimensionality of the normal
bundle presents the essential difficulties in such attempts.

\cite{mu1} studies the mean curvature flow of surfaces in a
Einstein four-manifold and  proves that a symplectic surface in a
K\"ahler-Einstein surface remains symplectic along the mean
curvature flow and the flow does not develope any type I
singularity. When the embient four-manifold $M$ is locally a
product of Riemann surfaces, there are two K\"ahler forms
$\omega'=\omega_1-\omega_2$ and $\omega''=\omega_1+\omega_2$ that
determine opposite orientations for $M$. We proved long time
existence and convergence of the mean curvature flow in \cite{mu1}
under the condition that the initial surface is symplectic with
respect to both $\omega'$ and $\omega''$. These results have been
generalized to arbitrary dimension and codimension in \cite{mu2}.
This article studies the case when the initial surface is
Lagrangian with respect to one K\"ahler form and symplectic with
respect to the other.

\vskip 10pt \noindent {\bf Theorem A} {\it Let $M$ be a compact
four manifold. If the universal covering of $M$ is any one of
$S^2\times S^2$, $\R^2\times \R^2$ or $\h^2 \times\h^2$ and
$\Sigma$ is a compact surface in $M$ that is Lagrangian with
respect to $\omega_1-\omega_2$ and symplectic with respect to
$\omega_1 +\omega_2$ , then the mean curvature flow of $\Sigma$
exists smoothly for all time. }
\vskip 10pt

It is proved by maximum principle in \S 2 that each slice
$\Sigma_t$ is again Lagrangian with respect to $\omega_1-\omega_2$
and symplectic with respect to $\omega_1+\omega_2$. In this
Lagrangian case, the condition of being symplectic with respect to
$\omega_1 +\omega_2$ is the same as saying $\Sigma$ is locally a
graph over both $\Sigma_1$ and $\Sigma_2$. This happens when
$\Sigma$ is the graph of a diffeomorphism between $\Sigma_1$ and
$\Sigma_2$. Recall a diffeomorphism $f:\Sigma_1 \rightarrow
\Sigma_2$ is called area preserving if $f^*\omega_2=\omega_1$. It
is not hard to see this is true if and only if the graph of $f$ in
$ M =\Sigma_1\times\Sigma_2 $ is an embedded Lagrangian surface
with respect to the symplectic structure
$\omega'=\omega_1-\omega_2$.

As for convergence at infinity, we prove the following general
subconvergence theorem.

\vskip 10pt \noindent{\bf Theorem B} {\it Let $\Sigma_1$ and
$\Sigma_2$ be compact Riemann surfaces with the same constant
curvature and $f:\Sigma_1 \mapsto \Sigma_2$ be an area preserving
diffeomorphism. As $t\rightarrow \infty$, a sequence of the mean
curvature flow of the graph of $f$ converges to a smooth minimal
Lagrangian graph. } \vskip 10pt

After this work was completed, the author was informed that K.
Smoczyk claims a proof to this theorem in the non-positive
curvature case assuming an extra angle condition.

The limit in this case is a "minimal map". This notion was
introduced by R. Schoen in \cite{sch}.

\begin{dfn}
A map $f:\Sigma_1 \rightarrow \Sigma_2$ is called a minimal map if
the graph is a minimal embedding in $M$.
\end{dfn}

Schoen also proved the existence and uniqueness of minimal
Lagrangian map when $\Sigma_1$ and $\Sigma_2$ are hyperbolic
surfaces. Theorem B gives a new proof of the existence part.

When $\Sigma_1=\Sigma_2$ and $f$ is homotopic to identity map, we
prove the following uniform convergence theorem.

\vskip 10pt \noindent {\bf Theorem C} {\it Let $\Sigma_1$ be a
compact Riemann surface of constant curvature and $f:\Sigma_1
\mapsto \Sigma_1$ be an area preserving diffeomorphism that is
homotopic to the identity map. Denote by $\Sigma$ the graph of $f$
in $\Sigma_1\times \Sigma_1$. The mean curvature flow of $\Sigma$
converges uniformly to a totally geodesic Lagrangian graph along
Lagrangian graphs . } \vskip 10pt

When $\Sigma_1$ is a sphere, we prove a stronger convergence in \S
3.

The result in \cite{mu1} has applications in the deformation of
maps between Riemann surfaces. It was proved that any map between
spheres with Jacobian less than one is deformed to a constant map
along the mean curvature flow of its graph. The results in this
article applies to the case when Jacobian is equal to one.

\vskip 10pt \noindent {\bf Corollary C} {\it Any area preserving
diffeomorphism $f :\Sigma_1\mapsto \Sigma_1$ that is homotopic to
the identity can be deformed to an isometry along area preserving
diffeomorphisms by the mean curvature flow. } \vskip 10pt

Since any diffeomorphism is isotopic to an area preserving
diffeomorphism, this gives a new proof of Smale's theorem
\cite{sma} that $O(3)$ is the deformation retract of the
diffeomorphism group of $S^2$. For a positive genus Riemann
surface, this implies the identity component of the diffeomorphism
group is contractible.

I would like to thank Professor R. Schoen
and Professor S.-T. Yau for their constant
encouragement and invaluable advice. I also
have benefitted greatly from the many discussion
that I have with Professor G. Huisken,
Professor L. Simon, and
Professor B. White.

\section{Long time existence}

Let $M$ be a smooth compact four manifold whose universal covering
is any of $S^2\times S^2$, $\R^2\times \R^2$ or $\h^2 \times\h^2$.
Equipping with the quotient metric, $M$ is in fact a locally
symmetric space. The standard K\"ahler forms on the factors extend
to two parallel forms $\omega_1$ and $\omega_2$ on $M$.
$\omega'=\omega_1-\omega_2$ and $\omega''=\omega_1+\omega_2$ are
two K\"ahler forms on $M$ that determine opposite orientations for
$M$ in the sense that $\omega'\wedge\omega'=
-\omega''\wedge\omega''$. The metric on $M$ is in particular
K\"ahler-Einstein with respect to either K\"ahler form and $Ric=c
\cdot g$. We shall fix the K\"ahler structure of $M$ to be
$\omega'$.

Let $F_0:\Sigma \rightarrow M$ be
a Lagrangian immersion of a compact surface
$\Sigma$.
We evolve $F_0$ in
the direction of its mean curvature.

\[\frac{dF}{dt}(x,t)=H(x,t)
\]
where $F:\Sigma \times [0,T) \rightarrow M$ is
a one parameter family of immersions $F_t(\cdot)
=F(\cdot, t)$ of $\Sigma$ and
$H(x,t)$ is the mean curvature vector of
$F_t(\Sigma)$ at $F_t(x)$. Whenever there is
no ambiguity, we shall write $F_t(\Sigma)
=\Sigma_t$.

Let $*$ be the Hodge operator on $\Sigma$,
then $*\omega_i$ is the Jacobian
of the projection from $\Sigma$ onto the $i$-th
factor of $M$.
\begin{dfn}
We say $\Sigma$ is a local graph with
respect to $\omega_i$ if $*\omega_i>0$ on
$\Sigma$.
\end{dfn}

The following proposition shows in particular
the condition of being the graph of an area
preserving diffeomorphism is preserved and thus
 the mean curvature flow does provide a deformation
for such diffeomorphisms.

\begin{pro}
Being a Lagrangian local graph in $M$ is preserved along the mean
curvature flow.
\end{pro}

\begin{proof}
The condition is equivalent to $\omega'=0$
and $\omega''>0$. By Proposition 4.1
in \cite{mu1} for any parallel K\"ahler
form $\omega$, $\eta=*\omega$ satisfies
the following equation,

\begin{equation}\label{eta1}
\begin{split}
\frac{d}{dt}\eta=\Delta\eta+\eta[(h_{31k}
-h_{42k})^2
+(h_{32k}+h_{41k})^2]+c \eta(1-\eta^2)
\end{split}
\end{equation}
\noindent
where $\{e_1, e_2, e_3,e_4\}$
is an orthonormal basis for $T_p M$ such that
$\{e_1, e_2\}$ forms an orthonormal basis
for $T\Sigma_t$. The basis is chosen
so that $d\mu(e_1, e_2)>0$ and
 $\omega^2(e_1, e_2, e_3, e_4)>0$ where $d\mu $ is
a fixed orientation on $\Sigma_t$.
Besides, $A(e_i,e_j)=h_{3ij}e_3+h_{4ij}e_4$
is the second fundamental form of $\Sigma_t$. From
this and the maximum principle
for parabolic equations we immediate see
that being a Lagrangian local graph is preserved.

\end{proof}

We fixed the complex structure $J'$
that corresponds to
$\omega'$ and choose a orthonormal basis
$\{e_1, e_2, e_3, e_4\}$ so that
$e_3=J'e_1$ and $e_4=J'e_2$. The orientation
given by this basis is in fact opposite
to the one given by $J'$ in the sense that
$(\omega')^2(e_1, e_2, e_3, e_4)<0$. Therefore
$(\omega'')^2(e_1, e_2, e_3, e_4)>0$ and
$\eta=*\omega''$ satisfies the equation
(\ref{eta1}).

 The normal bundle of a Lagrangian surface
 is canonically isomorphic to its tangent
 bundle by $J'$.
Through this isomorphism  the
 second fundamental form $A$
 and the mean curvature vector $H$ are
 associated with tensors
$B$ and $\sigma$.
$B$ is the symmetric three-tensor
defined by $B(X,Y, Z)=-<{\nabla}_X Y, J'(Z)>$
for $X, Y, Z \in T\Sigma$. $\sigma$ is
the one-form defined by
$\sigma(X)=<J'(X), H>$.  If
we denote $B(e_i, e_j, e_k)=B_{ijk}$ and
$\sigma(e_i)=\sigma_i$, then
$h_{3ij}=-B_{1ij}$ and $h_{4ij}=-B_{2ij}$. Therefore the term involving the second
fundamental form in (\ref{eta1}) can be calculated
as the following.

\[(h_{31k}-h_{42k})^2+(h_{32k}
+h_{41k})^2
=(B_{k11}-B_{k22})^2+4B_{k12}^2
=2|B|^2-|\sigma|^2
\]

Therefore $\eta=*\omega'' $ satisfies

\begin{equation}\label{eta}
\frac{d}{dt}\eta=\Delta \eta
+\eta[2|B|^2-|\sigma|^2]+c\eta(1-\eta^2)
\end{equation}

Now we proceed to prove the long time
existence theorem.

\vskip 10pt
\noindent
{\bf Theorem A}
{\it Let $M$ be a compact four manifold. If the universal covering of $M$ is any one
of $S^2\times S^2$,
$\R^2\times \R^2$ or $\h^2 \times\h^2$
and $\Sigma$ is a compact surface in $M$ that
is Lagrangian with respect to $\omega_1-\omega_2$
and symplectic with respect to $\omega_1
+\omega_2$ , then the mean curvature flow
of $\Sigma$ exists smoothly for all time.
}
\vskip 10pt

\begin{proof}

Notice that $0\leq \eta \leq 1$. By the
equation of $\eta$ and comparison
theorem for parabolic equations, we get

\[\eta(x,t)\geq \frac{\alpha e^{ct}}
{\sqrt{1+\alpha^2
e^{2ct}}}\]
where $\alpha>0$ satisfies $\frac{\alpha}
{\sqrt{1+\alpha^2}}
=\min_{\Sigma_0}\eta$. Therefore $\eta(x,t)$
converges uniformly to $1$ when $c=1$ and
is nondecreasing when $c=0$. In any case,
$\eta$ has a positive lower bound at any
finite time.

Since
\[|\sigma|^2=B_{111}^2
+2B_{111}B_{122}+B_{122}^2
+B_{211}^2+2B_{211}B_{222}+B_{222}^2\]

and

\[|B|^2=B_{111}^2+3B_{112}^2+3B_{122}^2
+B_{222}^2\]

It is easy to see $|\sigma|^2\leq \frac{4}{3}
|B|^2$, therefore

\[\frac{d}{dt}\eta \geq \Delta \eta +\frac{2}{3}\eta
|B|^2+c\eta(1-\eta^2)\]

We can proceed to prove regularity at
any finite time as in \cite{mu1}.
The idea is to prove the Gaussian density
\[ \lim_{t\rightarrow t_0}\int \rho_{y_0,t_0}
d\mu_t=1\]
for any point $y_0\in M$ and $t_0<\infty$,
where $\rho_{y_0, t_0}$ is the backward
heat kernel

\begin{equation}
\rho_{y_0, t_0}(y,t)=\frac{1}{(4\pi(t_0-t))}
\exp (\frac{-|y-y_0|^2}{4(t_0-t)})
\end{equation}
White's regularity theorem
\cite{w} would implies $(y_0,t_0)$ is a regular
point. As in \cite {mu1}, for any point
$(y_0,t_0) $ we can select sequences $t_i\rightarrow
t_0$ and $\lambda_i \rightarrow \infty$
such that the parabolic rescaling of
$\Sigma_{t_i}$ by $\lambda_i$ at $(y_0, t_0)$ converges to
a Lagrangian submanifold with $|B|=0$, or
a linear subspace. This implies the Gaussian
density at $(y_0, t_0)$ is $1$ and there is
 no singular point
at $(y_0, t_0)$.

\end{proof}

\section{Convergence at infinity-the sphere case}
In this section, we prove the convergence in the sphere case. The
key point is the uniform boundedness of the norm of the second
fundamental form. We accomplish this using the blow up analysis at
infinity. We already have a mean curvature flow $F:\Sigma\times
[0,\infty)\mapsto M$ that exists for all time. If
$\sup_{\Sigma_t}|A| $ is not bounded, then there exists a sequence
$\ct_k \rightarrow \infty$ such that
$\sup_{\Sigma_{\ct_k}}|A|\rightarrow \infty$. Choose $\cx_k\in
\Sigma_{\ct_k}$ such that $|A|(\cx_k, \ct_k)=\sup_{\Sigma_{\ct_k}}
|A|$. Fix a number $a$ less than the injective radius of $M$.
Because $M$ is compact, we may assume $\cx_k \rightarrow \cx \in
M$ and $d_M(\cx_k, \cx)< \frac{a}{2}$ by passing to a subsequence.

 Since
$M$ is locally a product, we can
choose a coordinate system on a neighborhood $U=U_1\times U_2$
of $\cx$ such that
each $\omega_i$ on $U_i$ is
the standard symplectic form. We
shall use this coordinate system
to identify $U$ with an open set $B$  in
$\R^4$. On $U$ there is the metric induced
from $M$ and on $B$ there is the flat
metric. However, being Lagrangian does
depend on any particular metric structure.

Let $S$ be the total space
of the mean curvature flow $F$  in $M\times[0,\infty)$,
take
\[S_k=S\cap (U\times[\ct_k -1, \ct_k +1])\]

For any $(x,t)\in S_k$ denote the parabolic
distance to the boundary of
$U\times [\ct_k -1, \ct_k +1] $ by
\[\delta_k (x,t)
=\min_{x_0\in \partial U, \,
t_0\in\{\ct_k -1, \ct_k +1\}}\{d_M (x, x_0),
\sqrt{|t-t_0|}\}\]

Denote
\[\alpha_k =\sup_{(x,t)\in S_k} \delta_k(x,t)
|A|(x,t)\]

Notice that $\alpha_k$ is a scaling invariant
quantity. Since
$\delta_k (\cx_k, \ct_k)|A|(\cx_k, \ct_k)
\geq \min\{\frac{a}{2}, 1\}
|A|(\cx_k, \ct_k)$, we have $\alpha_k \rightarrow
\infty$.
Now we consider $S_k$ as a smooth flow
in $B\subset \R^4$ with the  flat metric, let
\[\alpha'_k =\sup_{(x,t)\in S_k} \delta'_k(x,t)
|A'|(x,t)   \]
where
\[\delta'_k (x,t)
=\min_{x_0\in \partial B, \,
t_0\in\{\ct_k -1, \ct_k +1\}}\{|x- x_0|,
\sqrt{|t-t_0|}\}\]
is the parabolic distance in the flat metric
on $B$ and $|A'|(x,t)$ is the
second fundamental form of $\Sigma_t\cap B$ as
a submanifold in $\R^4$. Since the two
metrics are equivalent, $\alpha'_k\rightarrow
\infty$ too.

Now we take $(x_k, t_k)\in S_k$ such that

\[\delta'_k(x_k, t_k)|A'|(x_k, t_k)
\geq \frac{\alpha'_k}{2}\]

Let $\lambda_k=|A'|(x_k, t_k)$. Because
$\delta'_k (x_k, t_k)\leq \min \{a,1\}$,
$\lambda_k \rightarrow \infty$ too. For
any $S_k$, we consider it as a submanifold
in $B\times [\ct_k-1, \ct_k+1]
\subset \R^4\times \R$ and
take the parabolic rescaled flow
$\widetilde{S}_k=D_k S_k$ by $\lambda_k$.

\[
\begin{matrix}
&D_k :  &\R^4\times \R &\rightarrow
&\R^4\times \R\\
&&(x,t)&\rightarrow &(\lambda_k(x-x_k),
\lambda_k^2(t-t_k))
\end{matrix}
\]

Notice that $D_k(x_k, t_k)=(0,0)$. Let
$|\widetilde{A}_k|(y,s)$ denote
the second fundamental form of $ (\widetilde
{S}_k)_s$ at $y\in (\widetilde{S}_k)_s$, then
$|\widetilde{A}_k|(0,0)=1$.

Since  $\delta'_k|A'|$ is also a scaling invariant
quantity.

\[\widetilde{\delta}_k(0,0)|\widetilde{A}_k|
(0,0)=
\delta'_k(x_k, t_k)|A'|(x_k, t_k)
\geq \frac{\alpha'_k}{2}\]
where  $\widetilde{\delta}_k$ is the parabolic
distance to
the boundary of $\lambda_k(B-x_k)\times
[\lambda_k^2(\ct_k-1-t_k), \lambda_k(\ct_k+1
-t_k)]$.

From this we see $\widetilde{\delta}_k(0,0)
\rightarrow \infty$ as $k\rightarrow \infty$.
For any $(y,s)\in \widetilde{S}_k$, $\widetilde{\delta}_k(y,s)|\widetilde{A}_k|
(y,s)\leq \alpha'_k
\leq 2 \widetilde{\delta}_k(0,0)$, therefore

\[|\widetilde{A}_k|(y,s)\leq 2
\frac{\widetilde{\delta}_k (0,0)}
{\widetilde{\delta}_k (y,s) }
\leq 2\frac{\widetilde{\delta}_k (0,0)}
{\widetilde{\delta}_k (0,0)-\max\{|y|, \sqrt{|s|}\} }
\]

Take $(y, s)\in K$ for any compact set $K\subset \R^4\times \R$,
the above estimate shows $|\widetilde{A}_k|$ is uniformly bounded
for all $k$ on any compact set in space-time. Therefore
$\widetilde{S}_k \rightarrow \widetilde{S}_ \infty$ smoothly and
since $\widetilde{\delta}_k(0,0) \rightarrow \infty$,
$\widetilde{S}_\infty$ is defined on $(-\infty, \infty)$. We have
proved the following proposition.

\begin{pro}
If $\sup_{\Sigma_t}|A|$ is not bounded, then
there exists a blow-up flow $\widetilde{S}_
\infty \subset \R^4\times \R $ defined on
the whole $(-\infty, \infty)$ with
uniform bounded second fundamental form and
$|A|(0, 0)=1$.
\end{pro}

The main convergence theorem in the sphere case is the following.

\begin{thm}\label{sphere} Under the same assumption as in Theorem A. If $M$ has
positive curvature then the mean curvature flow of $\Sigma$
converges smoothly to a totally geodesic Lagrangian surface at
infinity.
\end{thm}

\begin{proof}
We already know the long time existence and
we are going to show the uniform boundedness of
the second fundamental form by contradiction.

By the equation of $\eta$ we have
$\eta(x,t)\geq 1-\epsilon_k$ for $(x,t)\in S_k$
and $\epsilon_k \rightarrow 0$. This continue
to hold for the corresponding $\widetilde{\eta}_k$
on $\widetilde{S}_k$. Therefore $\eta(x,t)
\equiv 1$  for $(x,t)\in \widetilde{S}_\infty
\subset \R^4\times\R$. In particular
the $t=0$ slice is a complete Lagrangian
graph with $\eta\equiv 1$.
If we write the graph as $(x,y, f(x,y), g(x,y))$,
then the Lagrangian condition implies

\[f_x g_y-g_x f_y=1\]
and $\eta=1$ is equivalent to
\[\frac{2}{\sqrt{1+f_x^2+f_y^2+g_x^2
+g_y^2+(f_x g_y-f_yg_x)^2}}=1\]
or
\[f_x^2+f_y^2+g_x^2
+g_y^2=2\]

These implies $h=f+\sqrt{-1} g$ is holomorphic
and $|\frac{\partial h}{\partial z}|=1$,
therefore $h$ is of the form $h=e^{\sqrt{-1}
\theta}z +C$, where $\theta$ and $C$ are
constants. The graph of $h$ has zero
second fundamental form and this
contradicts
with $|\widetilde{A}_\infty|(0,0)=1$.

Therefore $|A|^2$ is uniformly bounded in
space and time. By equation (7.7)
in \cite{mu1}, $|A|^2$ satisfies
the following equation.

\[\frac{d}{dt}|A|^2\leq
\Delta|A|^2-2|\nabla A|^2
+4|A|^4+K_1|A|^2+K_2
\]
Integrating this equation and we
see
\begin{equation}\label{bb}
\frac{d}{dt}\int_{\Sigma_t}
|A|^2 d\mu_t\leq C
\end{equation}
 Recall
$\frac{d}{dt}\eta  \geq \Delta\eta
+\frac{2}{3}\eta |A|^2$ and $\eta$ has
a positive lower bound, thus

\begin{equation}\label{cc}
\int_0^\infty \int_{\Sigma_t}
|A|^2 d\mu_t dt\leq \infty
\end{equation}

Equation (\ref{bb}) and (\ref{cc})
implies

\[\int_{\Sigma_t}|A|^2d\mu_t \rightarrow
0\]

By the small $\epsilon$ regularity theorem in
\cite{il}, $\sup_{\Sigma_t}|A|^2\rightarrow 0$ uniformly
as $t\rightarrow \infty$.

 Since the mean curvature
flow is a gradient flow and the metrics are
analytic, by the theorem of Simon \cite{si},
we get convergence at infinity. The flow converges to a minimal
Lagrangian submanifold with $\eta=1$. Since
$\eta=*(\omega_1+\omega_2)$, this implies
the limiting submanifold is holomorphic
with respect to the complex structure associated
with the K\"ahler form $\omega_1+\omega_2$.
The limiting map is both holomorphic and area
preserving and  thus an isometry.

\end{proof}

\section{Positive genus case}
The following theorem is a general subconvergence theorem for mean
curvature flow of surfaces. The proof is essentially contained in
that of Theorem 2 in \cite{si1}.

\begin{thm}\label{asym}
Let $F_t:\Sigma\times[0,\infty)\mapsto M$ be a smooth mean
curvature flow of an immersed compact oriented surface
$F_0(\Sigma)$ in a compact Riemannian manifold $M$. We assume
$F_0(\Sigma)$ represents a nontrivial homology class in $M$. Then
there exists a sequence $t_i \rightarrow \infty$ such that
$F_{t_i}$ converges to a $C^{1, \alpha}$ branched minimal
immersion $F:\widetilde{\Sigma}\mapsto M$ where
$\widetilde{\Sigma}$ is a compact oriented surface and $F$ is a
smooth immersion of $\widetilde{\Sigma}-B$ into $M$ for a finite
set $B\subset\widetilde{\Sigma}$.

\end{thm}

\begin{proof}
Since $\frac{d}{dt} \int d\mu_t=-\int |H|^2 d\mu_t$, we have
$\int_0^\infty (\int |H|^2 d\mu_t)dt <\infty$. Therefore there
exists a sequence $F_{t_i}$ such that $\int_{\Sigma_{t_i}} |H|^2
\rightarrow 0$. Now we apply Theorem 2 of \cite{si1} with the
functional $\frak{F}(\Sigma)=\frac{1}{2}\int_\Sigma |H|^2 d\mu
-\chi(\Sigma)$ where $\chi(\Sigma)$ is the Euler number of
$\Sigma$. By Gauss-Bonnet Theorem, this can be written as the form
of those functionals considered in \cite{si1}. $F_{t_i}$ form a
minimizing sequence of $\frak{F}$ in the space of immersions
smoothly homotopic to $F_0$. The homology class of $F_{t_i}$ is
nontrivial, so their diameters have positive lower bound. It
follows by Simon's theorem \cite{si1} that a limit exists as a
branched immersion.
\end{proof}

The convergence is in the sense of varifold and Housdorff distance
as discussed in  \cite{si2}. By \cite{si3} and Definition 2 of
\cite{si1}, the sequence $F_{t_i}$, while remains a minimizing
sequence for $\int |H|^2$, can be modified locally so that the
convergence is in the following sense. There exists a sequence
$\phi_{t_i}$ of diffeomorphisms of $\widetilde{\Sigma}-B$ onto
open subsets $U_{t_i}$ of $\Sigma$ such that
\begin{enumerate}
\item
$F_{t_i} \circ \phi_{t_i}$ converges to $F$ locally in the $C^2$
sense on $\widetilde{\Sigma}-B$.

\item
$F_{t_i}(M-U_{t_i})\subset \cup_{x\in B} B_{\epsilon_k}(F(x))$ for
some $\epsilon_k \downarrow 0$.
\end{enumerate}

By conclusion 2, it is not hard to see $F(\widetilde{\Sigma})$ is
in the same homology class as $F_0(\Sigma)$.  As was remark in
\cite{si1}, the surface $\widetilde{\Sigma}$ may have lower genus
than $\Sigma$ due to necks or handles pinching. Such pinching was
caused by the concentration of the limit measure
$|A|^2d\mu_{t_i}$. We shall prove Theorem B now.

\vskip .1in

\noindent {\bf Theorem B} {\it Let $\Sigma_1$ and $\Sigma_2$ be
compact Riemann surface with the same constant curvature and
$f:\Sigma_1 \mapsto \Sigma_2$ be an area preserving
diffeomorphism. As $t\rightarrow \infty$, a sequence of the mean
curvature flow of the graph of $f$ converges to a smooth minimal
Lagrangian graph.} \vskip .1in

\begin{proof}
The case when $\Sigma_1=\Sigma_2=S^2$ is already proved in Theorem
\ref{sphere}. We shall assume they are of positive genus now.
 By Theorem \ref{asym}, a subsequence
converges to a minimal Lagrangian immersion $F:
\widetilde{\Sigma}\mapsto \Sigma_1\times \Sigma_2$ which may
posses some branched points. We first show that indeed there is no
branched point. Since a minimal immersion is a conformal harmonic
map, the composite map $\pi_1\circ F:\widetilde{\Sigma} \mapsto
\Sigma_1$ is a harmonic map with respect to some smooth metric in
the same conformal calss as the pull back metric by $F$ . This is
now a degree one harmonic map since $F(\widetilde{\Sigma})$ is in
the same homology class as $\Sigma_{t_i}$ by conclusion 2 in the
remark right after Theorem \ref{asym}. Since the convergence
$F_{t_i}(\Sigma) \rightarrow F(\widetilde{\Sigma})$ is in the
varifold sense and $\pi_1\circ F_{t_i}$ has positive Jacobian,
$\pi_1\circ F$ has non-negative Jacobian. Use the proposition on
page 13 of \cite{sy} we can show the Jacobian of $\pi_1\circ F $
is positive everywhere and $\Sigma_\infty $ is the graph of a map
$f_\infty$. Therefore the sequence $\Sigma_{t_i}$ converges to a
smooth minimal Lagrangian graph. When $\Sigma_1$ and $\Sigma_2$
are both torus, the Gauss-Bonnet theorem shows $\int |A|^2
d\mu_{t_i} =\int |H|^2 d\mu_{t_i}\rightarrow 0$, so the limit is
totally geodesic.
\end{proof}

This gives a new proof of the existence theorem of minimal maps
between hyperbolic surfaces in Proposition 2.12 of \cite{sch}.

Indeed, Theorem B also holds when $\Sigma$ is locally a graph
which corresponds to the local condition
$\omega_1|_\Sigma=\omega_2|_\Sigma>0$.

\begin{thm}
Let $M=(\Sigma_1,\omega_1)\times(\Sigma_2, \omega_2)$, where
$\Sigma_1$ and $\Sigma_2$ are Riemann surfaces of the same
constant curvature. If $\Sigma$ is a compact Lagrangian surface
with respect to $\omega_1-\omega_2$ and is locally a graph over
$\Sigma_1$ and $\Sigma_2$, then the mean curvature flow of
$\Sigma$ exists for all time and a sequence converges to a smooth
minimal Lagrangian surface.
\end{thm}

\begin{proof}
The locally graphical condition implies
$\pi_1|_\Sigma:\Sigma\mapsto \Sigma_1 $ is a covering map, so
$2g-2=\deg(\pi_1|_\Sigma) (2g_1-2)$ where $g$ and $g_1$ is the
genus of $\Sigma$ and $\Sigma_1$ respectively. Now the limit
$\widetilde{\Sigma}$ has lower genus than $\Sigma$ and $\pi_1\circ
F $ is a branched harmonic immersion of degree
$\deg(\pi_1|_\Sigma)$ by the the same argument. This precludes the
possibility of branched point by the topological Riemann-Hurwitz
formula.
\end{proof}

We remark that the existence of minimal Lagrangian submanifold in
such homology class was first proved by Y.-I. Lee using
variational methods in \cite{lee}.

Next we prove Theorem C. First of all, we observe that when
$\Sigma_1=\Sigma_2$, the graph of the identity map is a totally
geodesic submanifold in the product space. When $f:\Sigma_1
\mapsto \Sigma_1$ is homotopic to the identity, we claim the limit
obtained in Theorem C is actually totally geodesic. When
$\Sigma_1$ is a torus, the Gauss-Bonnet theorem shows $\int |A|^2
d\mu_{t_i} =\int |H|^2 d\mu_{t_i}\rightarrow 0$, so the limit has
$|A|\equiv 0$. When $\Sigma_1$ is a hyperbolic surface, we apply
the uniqueness of minimal graph in each homotopy class in
\cite{sch} .

The next lemma should be well-known. We sketch the proof for
completeness.

\begin{lem}
Let $\Gamma$ be a totally geodesic submanifold in a Riemannian
manifold of non-positive sectional curvature and  $\rho(x)= d(x,
\Gamma)$ is the distance function to $\Gamma$. Then $\rho $ is a
convex function in a tubular neighborhood of $\Gamma$.
\end{lem}
\begin{proof}
Let $\alpha(s)$ be a smooth curve defined for
$-\epsilon<s<\epsilon$ and $\alpha (0)=x$. We need to show
$\frac{d^2}{ds^2}\rho(\alpha (s))\geq 0$. Join each point
$\alpha(s)$ to $\Gamma$ by a geodesic that realizes the distance
function. The lemma now follows from the second variation formula
of length (see e.g. page 20 in \cite{ce}) and the fact that
$\Gamma$ is totally geodesic.

\end{proof}

We slightly reformulate the statement of Theorem C and prove it in
the following. This formulation and Corollary C follows from the
correspondence between Lagrangian graph and area-preserving
diffeomorphism.

 \vskip .1in

\noindent {\bf Theorem C} {\it Let $\Sigma$ be a compact Riemann
surface of constant curvature and $f:\Sigma \mapsto \Sigma$ be an
area-preserving diffeomorphism  that is homotopic to the identity
map. The mean curvature flow $\Sigma_t$ of the graph of $f$ exists
for all time and each $\Sigma_t$ can be written as the graph of an
area-preserving diffeomorphism $f_t$. $f_t$ converges to the
identity map uniformly as $t \rightarrow \infty$.}

\vskip .1in

\begin{proof}
Let $\rho$ be the distance function to the diagonal in
$\Sigma\times\Sigma$. We calculate the parabolic equation of
$\rho(F(x,t))$ as the following.

\[\frac{d}{dt} \rho(F(x,t))=\nabla \rho \cdot H\]
where $\nabla \rho$ denotes the gradient of $\rho $ in
$\Sigma\times\Sigma$.

Split $\nabla \rho $ into the normal part and tangent part we get
$div_\Sigma \nabla \rho=\Delta_\Sigma \rho-\nabla \rho \cdot H$.

Therefore the equation is
\[ \frac{d}{dt} \rho =\Delta_\Sigma \rho-div_\Sigma \nabla \rho\]
$div_\Sigma \nabla \rho$ is the trace of the Hessian of $\rho$
restricted to $\Sigma$ and is always nonnegative by the convexity
of $\rho$ . By maximum principle, $\max_{x\in\Sigma} \rho(F(x,t))$
is non-increasing. Since we already have the convergence of a
subsequence $\Sigma_{t_i}$ in Hausdorff distance, this implies the
convergence of the whole flow $\Sigma_t$. We can write each
$\Sigma_t$ as $(x, f_t(x))$ for $x\in \Sigma$. Since $\rho((x,
f_t(x)) \rightarrow 0$ uniformly, this implies $f_t(x)$ converges
to $x$ uniformly.

\end{proof}

In fact, even when the domain and target of a diffeomorphism are
of different conformal structure, we can prove the following
isotopy theorem.

\begin{cor}
Every area preserving diffeomorphism of  Riemann surfaces is
isotopic through area preserving diffeomorphisms to a minimal
diffeomorphism.
\end{cor}
\begin{proof}
Given any area preserving diffeomorphism $h:\Sigma_1\mapsto
\Sigma_2$, we can compose it with another area-preserving minimal
diffeomorphism $g:\Sigma_2\mapsto \Sigma_1$ so that $g\circ h$ is
homotopic to the identity map on $\Sigma_1$. We can deform $g\circ
h$ by the mean curvature flow to get $f_t$ isotopic to the
identity map. Then $g^{-1}\circ f_t$ gives the desired isotopy of
$h$ to $g^{-1}$.

\end{proof}


\begin{thebibliography}{99}
\bibitem{ce} J. Cheeger and D. G. Ebin,
\textit{Comparison theorems in Riemannian geometry}, North-Holland
Mathematical Library, Vol. 9. North-Holland Publishing Co.,
Amsterdam-Oxford; American
 Elsevier Publishing Co., Inc., New York, 1975.
 viii+174 pp.



\bibitem{eh} K. Ecker and G. Huisken,
\textit{ Mean curvature evolution of
entire graphs.} Ann. of Math. (2) \textbf{130}
 (1989), no. 3, 453--471.



\bibitem{hu2} G. Huisken, \textit{Asymptotic
behavior for singularities of the mean curvature
flow}, J. Differential Geom. \textbf{31} (1990), no.
1, 285--299.


\bibitem{il} T. Ilmanen, \textit{Singularities
of mean curvature flow of surfaces },
preprint , 1997.

\bibitem{lee} Y.-I. Lee, \textit{Lagrangian minimal surfaces
in Kähler-Einstein surfaces of
negative scalar curvature}, Comm. Anal. Geom. 2 (1994), no. 4,
579--592.



\bibitem{sch} R. Schoen, \textit{
The role of harmonic mappings in rigidity and
 deformation problems}, Complex geometry
 (Osaka, 1990), 179--200, Lecture Notes in
 Pure and Appl. Math., 143, Dekker, New York,
 1993.
\bibitem{sy} R. Schoen and S.-T. Yau, \textit{Lectures on harmonic
maps}, Conference Proceedings and Lecture Notes in Geometry and
Topology, II. International Press, Cambridge, MA, 1997. vi+394 pp.
ISBN: 1-57146-002-0.

\bibitem{si1} L. Simon, \textit{Existence of Willmore surfaces}
Miniconference on geometry and partial differential equations
(Canberra, 1985),
 187--216, Proc. Centre Math. Anal. Austral. Nat. Univ., 10,
  Austral. Nat. Univ., Canberra, 1986.




\bibitem{si2} L. Simon, \textit{Existence of surfaces minimizing the Willmore functional}, Comm.
Anal. Geom.
 1 (1993), no. 2, 281--326.

\bibitem{si3} L. Simon, \textit{Personal communications}, April,
2001.


\bibitem{si} L. Simon, \textit{Asymptotics
 for a class of nonlinear evolution
 equations, with applications to geometric
  problems} Ann. of Math. (2) \textbf{118} (1983),
  no. 3, 525--571

\bibitem{sma} S. Smale, \textit{Diffeomorphisms
 of the $2$-sphere.} Proc. Amer. Math. Soc. 10
 1959, 621--626.

\bibitem{sm}  K. Smoczyk, \textit{
A canonical way to deform a Lagrangian
 submanifold.}, preprint, dg-ga/9605005



\bibitem{mu1}  M-T. Wang, \textit{Mean Curvature
Flow of surfaces in Einstein Four-Manifolds }, preprint , 2000.

\bibitem{mu2}  M-T. Wang, \textit{Stability
of graphic mean curvature flow in higher codimension }, preprint ,
2000.




\bibitem{w} B. White, \textit{A local
regularity theorem for classical mean curvature
flow}, preprint, 2000.



\end{thebibliography}
\end{document}